\documentclass[10pt]{article}
\usepackage[english]{babel}
\usepackage{amsfonts}
\usepackage[dvips]{graphicx}

\begin{document}

\title{\bf On a differential equation for a gas bubbles collapse mathematical model}

\date{June 2006}

\author{\bf Gianluca Argentini \\
\normalsize gianluca.argentini@gmail.com \\
\normalsize gianluca.argentini@riellogroup.com \\
\textit{Research \& Development Department}\\
\textit{Riello Burners}, 37048 San Pietro di Legnago (Verona), Italy}

\maketitle

\begin{abstract}
In this paper we present a mathematical model for estimate the collapse time of a gas bubble in a vane of a oil gerotor pump. This amount of time cannot be greater of the total time spent by the pump for filling and then emptying out a vane in a single revolution, otherwise there is a loss of lubrication between internal and external gears. We assume that oil is incompressible and viscous, the bubble has a spherical shape and it is not translating into the external fluid. The analytical treatment of the model shows that the Navier-Stokes equations for the velocity field of the oil can be reduced to a single non linear ordinary differential equation for the variation in time of the bubble radius. The collapse time estimated by a numerical resolution of this equation and the collapse time calculated from an analytical resolution of the linearized equation are substantially equal.\\

\noindent \textbf{keywords}: cavitation, Navier-Stokes equations, spherical coordinates, numerical and analytical resolution.
\end{abstract}

\section {The problem}

\noindent

\begin{figure}[ht]\label{pump}
	\begin{center}
	\includegraphics[width=6cm]{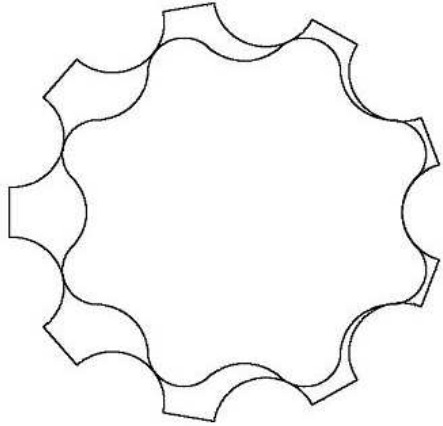}
	\caption{\small{\it Pignon and crown profiles in gerotor pump.}} 
	\end{center}
\end{figure}

We consider an oil gerotor-type pump for industrial applications, with an internal maximum pressure of about $15$ $bar$. Experimental results show that for a rotational velocity of about $2000$ $rpm$ some phenomena of cavitation can occur. In this work we don't consider a description of the arise of cavitation (see e.g. \cite{brennen}), but we want to estimate the time, which we call {\it collapse time} $t_c$, spent by a bubble of gas to reduce to zero its radius. In this way, we can deduce a maximum rotational velocity to avoid loss of lubrication due to the presence of bubbles in the geometrical zone where profiles of pignon and crown are coincident (see Figure \ref{pump}).

Assume that the gas (air or oil vapours) bubble has a spherical shape, with radius function of time: $R=R(t)$. Therefore the mathematical problem is to compute the time $t_c$ such that $R(t_c)=0$. The bubble is compressed by the external oil, which in general has a pressure $p$ greater than pressure $p_g$ of the internal gas. For simplicity, we assume that near the bubble the velocity field of the oil is radial, that is 

\begin{equation}\label{velocityField}
	\mathbf{v} = v(r,t)\mathbf{e}_r
\end{equation}

\noindent where $\mathbf{e}_r$ is the radial versor exiting from the centre of the bubble, $r$ is the radial distance of the oil material point from the centre, and $t$ is time variable (see Figure 2), and that the bubble is not translating into the fluid. Also, we don't consider external mass forces acting on the oil.

\begin{figure}[ht]
	\begin{center}
	\includegraphics[width=6cm]{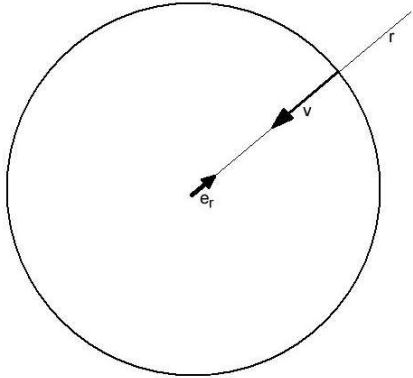}
	\caption{\small{\it Radial velocity field near the bubble.}} 
	\end{center}
\end{figure}

\noindent Our aim is to use Navier-Stokes equations for the oil flow near the bubble, so that we can use the velocity $v(R,t)$ at the bubble-boundary for computing the unknown function $R(t)$. Usual expression of these equations in a three dimensional cartesian system and without external forces is (see e.g. \cite{madani})

\begin{equation}\label{NavierStokes}
	\rho\left( \frac{\partial \mathbf{v}}{\partial t} + \mathbf{v} \cdot \nabla \mathbf{v} \right) = -\nabla p + \mu \Delta \mathbf{v}
\end{equation}

\noindent where $\rho$ is the density, $\mu$ the viscosity and $p$ the pressure of the fluid. For symmetry reasons due to the geometrical configuration, we'll see that this form is not the best for the treatment of the problem.

\section{Pressure of fluid near the bubble}

\noindent

Before using Navier-Stokes equations, we try to found an expression for the pressure $p$ of the fluid external to the bubble. For symmetry, we assume that near the bubble $p$ is a function of the radial distance $r$. Let $R_0$ the initial radius of the collapsing bubble.\\
Let $p_m$ the mean pressure, which we can assume to be constant, of the fluid in the pump. Then, for $r=R=0$, that is when the bubble is collapsed, $p$ is equal to $p_m$, while for $r=R=R_0$ $p$ can be considered as the difference between $p_m$ and the pressure $p_{g0}$ due to the gas contained on a sphere of radius $R_0$. From the {\it perfect gas law} we have

\begin{equation}\label{pgl}
	p_{g0} = \frac{3n_0\mathbf{R}T}{4 \pi R_0^3}
\end{equation}

\noindent where $n_0$ is the number of gas moles, $\mathbf{R}$ is the gas universal constant and $T$ the gas temperature. If we assume that $p$ is linear on $r$, the simplest expression for $p$ to satisfy the previous conditions is

\begin{equation}\label{pExpressTemp}
	p = cRr + \frac{R_0 - R}{R_0}p_m
\end{equation}

\noindent where $c$ is a constant which can be computed imposing the condition $p = p_m -p_{g0}$ when $r=R=R_0$. The result is

\begin{equation}\label{pExpress}
	p = \left(\frac{p_m}{R_0^2} - \frac{3n_0\mathbf{R}T}{4 \pi R_0^5} \right)Rr + \frac{R_0 - R}{R_0}p_m
\end{equation}

\noindent Note that, as expected, $p$ is a radial symmetry function.

\section{Use of spherical coordinates system}

\noindent The physical nature of the problem and its mathematical modelization induce us to substitute the cartesian coordinates with a spherical coordinates system with origin at the centre of the bubble. These coordinates are the radial distance $r$ from the origin, the azimuthal angle (longitude) $\alpha$, with $0 \leq \alpha \leq 2\pi$, and the elevation from the plane of the equator (latitude) $\theta$, with $0 \leq \theta \leq \pi$ (see {\cite{arfken}).\\
We have to translate the Navier-Stokes equations (\ref{NavierStokes}) from cartesian to new spherical coordinates in the case of radial velocity field $\mathbf{v}(\mathbf{x},t)=v(r,t)\mathbf{e}_r$ and radial pressure $p=p(r,t)$. Note that $r^2 = x^2 + y^2 + z^2$, therefore

\begin{equation}\label{radialCoordGradient}
	\partial_x r = \frac{x}{r}, \hspace{0.1cm} \partial_y r = \frac{y}{r}, \hspace{0.1cm} \partial_z r = \frac{z}{r}
\end{equation}

\noindent Also we have 

\begin{equation}
	\mathbf{e}_r=\frac{\mathbf{x}}{r}
\end{equation}

\noindent Then, denoted by $v_i$ the $i$-th cartesian component of the velocity $\mathbf{v}$, from the chain rule the following formula for the $k$-th component of the gradient of $v_i$ holds: 

\begin{equation}\label{gradientSpherical}
	\partial_k v_i = \frac{\partial v}{\partial r} \frac{x_k^2}{r^2} + v\left[ \frac{r^2\delta_{ki}-x_kx_i}{r^3} \right]
\end{equation}

\noindent From these identities, the cartesian divergence $\nabla \cdot \mathbf{v}$ has the following shape in the new coordinates:

\begin{equation}\label{divergenceSpherical}
	\nabla \cdot \mathbf{v} = \sum_{i=1}^3 \partial_i v_i = \frac{\partial v}{\partial r} + \frac{2}{r}v
\end{equation}

\noindent Now we apply the hypothesis of incompressibility for the fluid; from $\nabla \cdot \mathbf{v} = 0$ follows

\begin{equation}\label{incompressTemp}
	\frac{\partial v}{\partial r} = - \frac{2}{r}v
\end{equation}

\noindent which can be resolved by variables separation, so that

\begin{equation}\label{velocityRadialTemp}
	v(r,t) = \frac{f(t)}{r^2}
\end{equation}

\noindent with $f$ arbitrary function of $t$. If we assume that the bubble boundary has velocity $\dot{R}(t)$ equals to the external fluid velocity field when $r=R(t)$ (physically, this corresponds to the fact that fluid doesn't enter into the bubble deforming its geometrical shape), we have

\begin{equation}\label{ftTemp}
	\dot{R}(t) = v(R(t),t)=\frac{f(t)}{R^2(t)}
\end{equation}

\noindent and therefore $f(t)=R^2(t)\dot{R}(t)$. Then, from (\ref{velocityRadialTemp}), the expression for the velocity field is

\begin{equation}\label{velocityRadial}
	\mathbf{v}(r,t) = \frac{R^2(t)\dot{R}(t)}{r^2}\mathbf{e}_r
\end{equation}
	
\noindent Now we can write the partial derivative respect time variable, necessary in Navier-Stokes equations:

\begin{equation}\label{timeDerivativeVelocity}
	\frac{\partial \mathbf{v}}{\partial t} = \frac{1}{r^2}\left(R^2\ddot{R}+2R\dot{R}^2\right)\mathbf{e}_r
\end{equation}

\noindent Then we translate in spherical coordinates the cartesian expression $\mathbf{v} \cdot \nabla \mathbf{v} = \left(\sum_{k=1}^3 v_k \partial_k v_i\right)_{1\leq i \leq3}$. From identities (\ref{gradientSpherical}) and from $v_k = v(r,t)(\mathbf{e}_r)_k = v \frac{x_k}{r}$, we obtain

\begin{equation}
	v_k \partial_k v_i = \frac{\partial v}{\partial r} v \frac{x_k^3}{r^3} + v^2 \left[ \frac{r^2 x_k \delta_{ki}-x_k^2 x_i}{r^4} \right]
\end{equation}

\noindent Summing on $k$ and using $x^2+y^2+z^2=r^2$, $\mathbf{e}_r=\frac{\mathbf{x}}{r}$, the following expression holds:

\begin{equation}\label{nonLinearTermSphericalTemp}
	\mathbf{v} \cdot \nabla \mathbf{v} = \frac{\partial v}{\partial r} v \left(\frac{x}{r},\frac{y}{r},\frac{z}{r}\right) = \frac{\partial v}{\partial r} v \mathbf{e}_r
\end{equation}

\noindent From (\ref{velocityRadial}) and (\ref{incompressTemp}) we can deduce the new form of the convective term in fluid dynamics equations:

\begin{equation}\label{nonLinearTermSpherical}
	\mathbf{v} \cdot \nabla \mathbf{v} = \left( -\frac{2}{r^5}R^4(t)\dot{R}^2(t) \right) \mathbf{e}_r
\end{equation}

\noindent Now, using standard computations (or see (\cite{feynman}) for a calculation in the case of spherical symmetry), from (\ref{gradientSpherical}) and (\ref{incompressTemp}) we can deduce the expression for the laplacian of $\mathbf{v}$ in spherical coordinates:

\begin{equation}\label{laplacianSpherical}
	\Delta \mathbf{v} = \sum_{k=1}^3 \partial^2_{kk} \mathbf{v} = \left(\frac{\partial^2 v}{\partial r^2} + \frac{2}{r}\frac{\partial v}{\partial r}\right)\mathbf{e}_r = \frac{2}{r^2}v\mathbf{e}_r = \frac{2}{r^4}R^2(t)\dot{R}(t)\mathbf{e}_r
\end{equation}

\noindent Finally, we derive the form for the gradient of the pressure $p$:

\begin{equation}\label{gradientPressureSphericalTemp}
	(\nabla p)_i = \frac{\partial p}{\partial r}\frac{x_i}{r}= \left(\frac{p_m}{R_0^2} - \frac{3n_0\mathbf{R}T}{4 \pi R_0^5} \right)R \frac{x_i}{r}
\end{equation}

\noindent so that

\begin{equation}\label{gradientPressureSpherical}
	\nabla p = \left(\frac{p_m}{R_0^2} - \frac{3n_0\mathbf{R}T}{4 \pi R_0^5} \right)R \mathbf{e}_r
\end{equation}

\section{An ordinary differential equation}

\noindent 

In this section we derive an ordinary differential equation which translate Navier-Stokes system into a new form using radial symmetry and the spherical coordinates.\\
We substitute the differential expressions in original fluid dynamics equations with their new forms computed in the previous section. For simplicity, let $a=\frac{p_m}{R_0^2} - \frac{3n_0\mathbf{R}T}{4 \pi R_0^5}$, with $a$ constant. Therefore, using (\ref{timeDerivativeVelocity}), (\ref{nonLinearTermSpherical}), (\ref{laplacianSpherical}) and (\ref{gradientPressureSpherical}) in equations (\ref{NavierStokes}), we can write

\begin{equation}
	\rho \left( \frac{R^2\ddot{R}+2R\dot{R}^2}{r^2} - \frac{2}{r^5}R^4\dot{R}^2 \right) = -aR + \mu \frac{2}{r^4}R^2\dot{R}
\end{equation}

\noindent For $r=R$ this equation describe the evolution of the bubble radius; simplifying similar terms and multiplying by $R^2$, we obtain the following equation for the unknown function $R(t)$

\begin{equation}\label{odeNS}
	\rho\ddot{R}R^2 - 2\mu\dot{R} + aR^3 = 0
\end{equation}

\noindent Note that the leading differential term $\ddot{R}R^2$ is the same than that of {\it Rayleigh-Plesset} equation of bubble dynamics in cavitation (see \cite{brennen}). The previous {\it non linear} second order equation seems not analytically integrable with classical methods; probably an hard calculation using {\it Lie symmetries} technique might be useful for informations about expression of an exact solution (see e.g. \cite{stephani}). But in the next sections we deduce results of industrial interest using numerical and approximate analytical resolution.\\
We want to solve an initial value problem for (\ref{odeNS}) with conditions

\begin{equation}
	R(0)=R_0, \hspace{0.1cm} \dot{R}(0)=0
\end{equation}

\noindent The condition on derivative seems quite realistic. Note that a not null constant function $R(t)=k$ $\forall t$ is not a solution of the equation, as expected for physical and technical reasons.

\section{Physical considerations}

\noindent

For a computation of the coefficient $a$ in equation (\ref{odeNS}) the main step is the estimate of the number of moles. As first, but quite realistic, approximation we consider that the gas contained in the bubble is atmospheric air. Let $R_0=0.05$ $cm$ the initial radius of the bubble, a value compatible with the dimensions of a body pump of industrial interest. Using the value $\rho_{air}=0.01177$ $g/cm^3$ as air density and the value $W_{air}=28.97$ $g/mol$ as air molecular weight, we calculate the air mass $M_{air}$ contained in a sphere of radius $R_0$, so we can know the number of moles $n_0=\frac{M_{air}}{W_{air}}$ (see e.g. \cite{atkins} for details). The resulting value of $a$ is about $10^7$ $dyne/cm^4$.
As fluid we use light oil with $\rho=8.2$ $g/cm^3$ and $\mu=0.0287$ $St$ (at standard temperature of $300^{\circ}K$). As mean pressure of fluid into vanes we consider the value $p_m=10^7$ $dyne/cm^2$ $\approx$ $10$ $bar$, that is the value usually reached at the halfway during the angular rotation of a single vane pump.

\section{Numerical resolution}

\noindent

Consider now the numerical resolution of our differential problem. We have used the software {\it Mathematica}, ver. 5.2, and its built-in function {\ttfamily NDSolve}, which applies adaptive finite differences schemes like Adams and Backwards Differentiation (see \cite{shampine} for informations about these methods) with high degree of accuracy (see \cite{coombes} for details). The result is shown in Figure 3.\\
The collapse time $t_c$ is about $0.00141$ $sec$. At a frequence of $2000$ $rpm$, during an interval of time equal to $t_c$ the pignon rotates of about $17^{\circ}$, therefore if bubble is produced near the geometrical zone where pignon and crown have conciding profiles (right zone of Figure \ref{pump}), the probability of loss of lubrication and physical contact between rotating parts can be not null.

\begin{figure}[ht]
	\begin{center}
	\includegraphics[width=6cm]{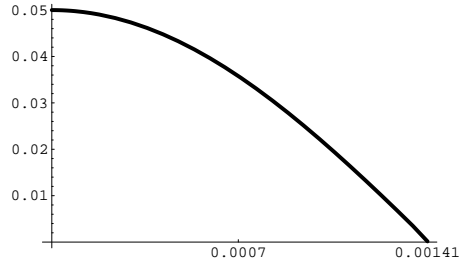}
	\caption{\small{\it Profile for numerical solution; variable in horizontal axis is time in $sec$.}} 
	\end{center}
\end{figure}

\section{Analytical resolution}

\noindent

In this section we approach an analytical approximated resolution for equation (\ref{odeNS}). From previous numerical resolution the key information is that $t_c << 1$, therefore an approximation based on some type of expansion respect the time variable and centered at $t=0$ is well justified. For this aim we consider an application of {\it Implicit Function} Theorem to the case of symbolic resolution of differential equations (see \cite{argentini} for details).\\
From hypothesis we have $R(0)=R_0$ and $\dot{R}(0)=0$, so from (\ref{odeNS}) follows that

\begin{equation}
	\ddot{R}(0)=-\frac{aR_0}{\rho}
\end{equation}

\noindent Consider the polynomial function $F(t,R,\dot{R},\ddot{R})=\rho\ddot{R}R^2-2\mu\dot{R}+aR^3$, where $R$, $\dot{R}$ and $\ddot{R}$ are to be considered independent variables and not functions. From previous formula, we have 

\begin{equation}
	F(0,R_0,0,-\frac{aR_0}{\rho})=0
\end{equation}

\noindent Then we compute all the partial derivatives of $F$ at point $(0,R_{\hspace{0.05cm}0},0,-\frac{aR_0}{\rho})$:

\begin{equation}
	\partial_t F = 0, \hspace{0.1cm} \partial_R F = aR_0^2, \hspace{0.1cm} \partial_{\dot{R}} F = -2\mu, \hspace{0.1cm} \partial_{\ddot{R}} F = \rho R_0^2
\end{equation}

\noindent As $\partial_{\ddot{R}} F \neq 0$, from implicit function Theorem (see e.g. \cite{courant}) we can write the following approximated but explicit formula for the variable $\ddot{R}$

\begin{equation}
	\ddot{R} = -\frac{aR_0}{\rho} - \frac{a}{\rho}(R-R_0) + \frac{2\mu}{\rho R_0^2}(\dot{R}-0)
\end{equation}

\noindent from which we obtain a linear constant coefficients ordinary differential equation for the function $R(t)$:

\begin{equation}\label{odeNSlinear}
	\rho R_0^2\ddot{R} - 2\mu\dot{R} + aR_0^2R = 0
\end{equation}

\noindent The analytical solution of the initial value problem can be found with the classical method of the characteristic equation, and its shape is quite complicated:

\begin{equation}
	\frac{e^{\frac{\left(\mu -\sqrt{\mu ^2-aR_0^3 \rho }\right) t}{R_0^2 \rho }}R_0 \left[\sqrt{\mu ^2-aR_0^3 \rho } \left(1+e^{\frac{2
   \sqrt{\mu ^2-aR_0^3 \rho } t}{R_0^2 \rho }}\right)-e^{\frac{2 \sqrt{\mu ^2-aR_0^3 \rho } t}{R_0^2 \rho }} \mu +\mu \right]}{2 \sqrt{\mu
   ^2-aR_0^3 \rho }}
\end{equation}

\noindent But numerical experiments testify that its second order Taylor expansion centered at $t=0$ gives a very accurate approximation, so that we can establish that the shape of the approximated analytical solution of the differential problem of bubble collapse is

\begin{equation}\label{analyticalSolution}
	R(t) = R(0)+\dot{R}(0)(t-0)+\frac{1}{2}\ddot{R}(0)(t-0)^2 = \frac{aR_0}{2\rho}t^2 + R_0
\end{equation}

\noindent The collapse time $t_c$ in this case is given by resolution of the algebraic equation $R(t_c)=0$:

\begin{equation}\label{analyticalCollapseTime}
	t_c = \sqrt{\frac{2\rho}{a}}
\end{equation}

\noindent Substituting the values of $\rho$ and $a$, we found $t_c = 0.00128$ $sec$, which matches the previous value estimated by numerical resolution until millisecond. Note that viscosity $\mu$, in this approximation, is not relevant for computing collapse time. If a Taylor expansion of third order is used for the analytical solution of (\ref{odeNSlinear}), the coefficient of $t^3$ is

\begin{equation}\label{cubicTerm}
	\frac{1}{6}R^{(3)}(0) = a_3 = - \frac{a\mu}{3\rho^2R_0}
\end{equation}

\noindent and the numerical contribute of the cubic term to reduce $R(t)$ is very small, $|a_{\hspace{0.05cm}3}t^3| < 10^{-4}$ $cm$.\\
In Figure 4 we compare the two profiles of $R(t)$ obtained by numerical and analytical methods.

\begin{figure}[ht]
	\begin{center}
	\includegraphics[width=6cm]{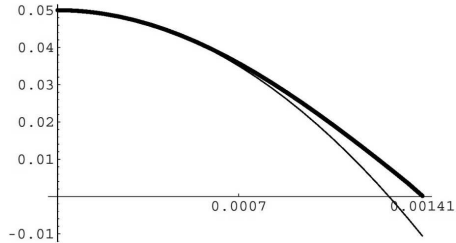}
	\caption{\small{\it Light line is the analytical profile.}} 
	\end{center}
\end{figure}

\section{Conclusions}

\noindent

The analytical resolution gives a formula for estimate the collapse time, and in particular from (\ref{analyticalCollapseTime}) and (\ref{cubicTerm}), we can deduce that a fluid with a lower density, or with a greater viscosity, or a pump with a greater mean pressure in vanes can give benefits for a reduction of time permanence of bubbles. The approximated calculations in analytical case give results matching very well those of numerical resolution. But in the second case the previous technical informations about density and pressure can be obtained only by an amount of computational simulations.\\
The mathematical model can be improved, in particular assuming for the bubble a geometrical shape more complicated than the spherical one. In this case the function $f(t)$ of equations (\ref{velocityRadialTemp}) and (\ref{ftTemp}) should have a more complicated dependence on $R(t)$, therefore the differential equation could become more difficult to resolve.\\
The non linearity of the differential equation (\ref{odeNS}) induces to consider more subtle techniques for searching a non approximated analytical resolution. Probably the method of Lie symmetries generating a two dimensional Lie algebra could transform the original equation to an integrable form.

\end{document}